\documentclass[12pt]{article}

\usepackage{amscd,amssymb}
\usepackage{amsfonts}
\usepackage{theorem,url}
\usepackage[leqno]{amsmath}
\usepackage{enumerate}
\usepackage{indentfirst, color}
\usepackage{hyperref}

\newtheorem{thm}{Theorem}
\newtheorem{defi}{Definition}
\newtheorem{rem}{Remark}
\newtheorem{ex}{Example}
\newtheorem{prop}{Proposition}

{\noindent  \textbf{#2}} %
{\hfill $\diamond$}

   \newcommand{\dd}{
      \mathop{}\mathopen{}\mathrm{d}
   }

\newcommand{\Prb}{\mathbb{P}}
\newcommand{\E}{\mathbb{E}}

\newcommand*{\cA}{\mathcal{A}}

\newcommand*{\cD}{\mathcal{D}}

\newcommand*{\cF}{\mathcal{F}}

\newcommand*{\cK}{\mathcal{K}}

\newcommand*{\cM}{\mathcal{M}}

\newcommand*{\cR}{\mathcal{R}}
\newcommand*{\cS}{\mathcal{S}}

\newcommand*{\cX}{\mathcal{X}}

\newcommand{\bB}{\mathbb{B}}
\newcommand{\bC}{\mathbb{C}}
\newcommand{\bD}{\mathbb{D}}
\newcommand{\bE}{\mathbb{E}}
\newcommand{\bF}{\mathbb{F}}

\newcommand{\bH}{\mathbb{H}}

\newcommand{\bL}{\mathbb{L}}
\newcommand{\bM}{\mathbb{M}}

\newcommand{\bP}{\mathbb{P}}

\newcommand{\bR}{\mathbb{R}}
\newcommand{\bS}{\mathbb{S}}

\newcommand{\eps}{\varepsilon}

\newcommand{\ind}{\mathbf{1}}

\newcommand{\tr}{\text{{\rm Tr}}}

\newcommand{\supp}{\text{supp}}
\newcommand{\cad}{c\`adl\`ag }

\title{Limit behaviour of the minimal solution of a BSDE in the non Markovian setting.}

\author{Dmytro Marushkevych\thanks{Laboratoire Manceau de Math\'ematiques, Le Mans Universit\'e, Avenue Olivier Messiaen, 72085 Le Mans cedex 9, France.\hfill\break
e-mail: {\tt Dmytro.Marushkevych.Etu@univ-lemans.fr, Alexandre.Popier@univ-lemans.fr}} \ and 
Alexandre Popier\footnotemark[1]
}

\begin{document}

\maketitle

\begin{abstract} 
We use the functional It\^o calculus to prove that the solution of a BSDE with singular terminal condition is continuous at the terminal time. Hence we extend known results for a non-Markovian terminal condition.
\end{abstract}

\noindent {\bf AMS class:} 60G99, 60H99.

\section*{Introduction}

In this paper we consider a filtered probability space $(\Omega,\cF, \bF,\Prb)$ with a complete and right-continuous filtration $\bF=\{\cF_t, \ t \geq 0\}$. We assume that this space supports a Brownian motion $W$. We consider the following BSDE:
\begin{equation}\label{eq:BSDE_jumps}
Y(t)=\xi + \int_t^T f(s,Y(s),Z(s)) ds - \int_t^T Z(s) dW(s) - \int_t^T dM(s)
\end{equation}
where $f$ is the generator and $\xi$ is the terminal condition. The solution is the triplet $(Y,Z,M)$. Since no particular assumption is made on the underlying filtration, there is the additional martingale part $M$ orthogonal to $W$. It is already established that such a BSDE has a unique solution when the terminal condition $\xi$ belongs to $L^p(\Omega,\mathcal{F}_T,\mathbb{P})$, $p >1$ (see among others \cite{delo:13} or \cite{krus:popi:16}). 

When the terminal condition $\xi$ satisfies
\begin{equation} \label{eq:sing_term_cond}
\Prb(\xi = +\infty)  >0
\end{equation}
we called the BSDE {\it singular}. This singular case has been studied in \cite{popi:06} when the filtration is generated by the Brownian motion (no additional noise, i.e. $M=0$) and for the particular generator $f(t,y,z)=f(y)=-y|y|^q$. Extension have been studied in \cite{krus:popi:16} and \cite{popi:16}. Recently singular BSDE were used to solve a particular stochastic control problem with application to portfolio management (see \cite{anki:jean:krus:13}, \cite{grae:hors:qiu:13}, \cite{krus:popi:16b}). In this framework, the generator does not depend on $z$ and has the following form:
\begin{equation} \label{eq:example_control}
f(t,y,z) = - \frac{y|y|^{q}}{q \alpha(t)^{q}} + \gamma(t).
\end{equation}
where $\alpha$, $\beta$ and $\gamma$ are positive processes. The minimal solution $(Y,Z,M)$ (provided it exists) gives the value function of the following control problem: minimize\footnote{with the convention $0.\infty = 0$.}
\begin{equation}\label{eq:control_pb_intro}
\E \left[  \int_t^T \left( \alpha(s) |\eta(s)|^p + \gamma(s) |X(s)|^p \right) ds  + \xi  |X(T)|^p \bigg| \cF_t \right]
\end{equation}
over all progressively measurable processes $X$ that satisfy the dynamics
\begin{equation*}
X(s) =x +\int_t^s \eta(u) du 
\end{equation*}
and the terminal state constraint 
\begin{equation*}
X(T) \ind_{\xi=\infty}= 0.
\end{equation*}
$p$ is the H\"older conjugate of $1+q$. For the financial point of view, the set $\{\xi=+\infty\}$ is a specification of a set of market scenarios where liquidation is mandatory. The value function is equal to $|x|^p Y_t$ and the optimal state process $X^*$ can be computed directly with $Y$. Note that the martingale part of the solution $(Z,M)$ is not employed in the computation of the optimal state process. Thus the control problem can be completely solved provided the BSDE has a minimal solution (see Section 2 and Theorem 4 in \cite{krus:popi:16b} for more details on the control problem).

In \cite{krus:popi:16b}, under some technical sufficient assumptions on $f$ (Conditions \textbf{(A)} below), it is proved that the BSDE \eqref{eq:BSDE_jumps} with singular terminal condition \eqref{eq:sing_term_cond} has a minimal super-solution $(Y,Z,M)$ such that a.s.
\begin{equation} \label{eq:term_cond_super_sol}
\liminf_{t\to T} Y(t) \geq \xi.
\end{equation}
The main requirement is that $f$ decreases w.r.t. $y$ at least as a polynomial function (almost like $-y^{1+q}$, $q>0$), when $y$ is large. The main difficulty is to obtain some a priori estimate, which states that $Y_t$ is bounded from above for any $t < T$ by a finite process (Inequality \eqref{eq:a_priori_estimate}). One can construct the solution $(Y,Z,M)$ without Condition \eqref{eq:term_cond_super_sol} if the filtration $\bF$ is complete and right-continuous. 

In the classical setting ($\xi \in L^p(\Omega)$), $Y$ has a limit as $t$ increases to $T$ since the process is solution of the BSDE \eqref{eq:BSDE_jumps} and thus is c\`adl\`ag\footnote{French acronym for right continuous with left limits.}. Moreover this limit is equal to $\xi$ a.s. if the filtration $\bF$ is left-continuous at time $T$. Indeed we need to avoid a jump at time $T$ of the orthogonal martingale $M$ (see the later discussion in Section \ref{sect:setting}). Hence in the singular case the behaviour \eqref{eq:term_cond_super_sol} of the super-solution $Y$ at time $T$ is obtained under this additional requirement on the filtration $\bF$. For the related control problem \eqref{eq:control_pb_intro}, this weak behaviour \eqref{eq:term_cond_super_sol} at time $T$ of the minimal process $Y$ is sufficient to obtain the optimal control and the value function (see \cite{krus:popi:16b}). Nevertheless two natural questions arise here:
\begin{enumerate}
\item Does the limit exist ? 
\item Can the inequality \eqref{eq:term_cond_super_sol} be an equality if the filtration is left-continuous at time $T$ ?
\end{enumerate}
Despite the very theoretical aspect of these questions, there are several applications. From the financial point of view it means that the optimal liquidation portfolio does not super hedge the penalty cost $\xi$. And in \cite{bank:voss:18}, a positive answer to these questions is a condition for solving the optimal targeting problem.

\subsubsection*{Related literature}

As far as we know, there are only three works on this topic: \cite{popi:06}, \cite{popi:16} and \cite{krus:popi:seze:18}. In \cite{popi:06} we were able to prove this in the Brownian setting, that is when the filtration $\bF$ is generated by the Brownian motion $W$ and if $f(t,y,z) = -y|y|^q$. We proved that the limit always exists (see \cite[Proposition 9]{popi:06}) and to obtain the equality, we supposed that $\xi = g(X(T))$ where $X$ is the solution of a forward SDE. We distinguished two cases:
\begin{itemize}
\item When $q>2$ without additional conditions since we have a suitable control of $Z$.
\item When $q\leq 2$ but with Malliavin calculus: roughly speaking $Z$ is the Malliavin derivative of $Y$ and we use the integration by parts to remove $Z$. 
\end{itemize}

In \cite{popi:16}, we deal with a generator satisfying Condition {\bf (A)} (see below). The filtration $\bF$ should be left continuous at time $T$ (to avoid thin time case, see the discussion in \cite[Section 2.2]{popi:16}). 
\begin{enumerate}
\item The existence of a limit at time $T$ is proved under a structural condition on the generator $f$ (\cite[Theorem 3.1]{popi:16}). Roughly speaking we prove that $Y$ is a non linear continuous transform of a non negative supermartingale. Relaxing the condition on $f$ is not the aim of this paper. 
\item In \cite{popi:16}, we also extend the result concerning the second question, again when $\xi = g(X(T))$. This setting is called {\it half-Markovian} since we do not require any similar condition on $f$. 
\end{enumerate}

The paper \cite{krus:popi:seze:18} was the first attempt to prove the equality when $\xi$ is not given by $g(X(T))$ in the Brownian setting. Indeed $\xi$ was assumed to be equal to $\infty \mathbf 1_{B(m,r)^c}$ or $ \infty \mathbf 1_{B(m,r)}$, where $B(m,r)$ is the ball in the space $C([0,T])$ of the underlying Brownian motion centered at the constant function $m$ and radius $r$. Our proof was based on the exit time of the Brownian motion and the derivation and solution of a related heat equation with a singular and discontinuous Dirichlet boundary condition. Let us remark that the considered functionals are not continuous, in the sense of \cite[Definition 2.3]{cont:four:13}.

\subsubsection*{Contributions and decomposition of the paper}

Our goal is to give another class of non Markovian terminal values $\xi$ such that \eqref{eq:term_cond_super_sol} becomes:
$$\liminf_{t\to T} Y(t) = \xi.$$
This class is constructed using the functional It\^o calculus developed by \cite{cont:16,cont:four:10,cont:four:13,dupi:09} (see \cite{cont:16} for an overview and the references therein). Roughly speaking, $\xi$ is a smooth functional $F$ of the paths of a continuous diffusion process $X$ (solution of the SDE \eqref{eq:SDE}) and of its bracket $[X]$, satisfying some integrability  assumption (see Condition {\bf (C)}). As presented in the subsection \ref{ssect:examples}, our condition includes the Markovian case studied in \cite{popi:16}, but also the integral of $X$ w.r.t. $t$ or some approximation of the process $X$. 

\bigskip
The paper is decomposed as follows. In Section \ref{sect:setting}, we recall known results concerning BSDEs with singular terminal condition and the functional It\^o calculus, in particular the definition of smooth functionals. In Section \ref{sect:continuity_Y_Markov}, we give the setting of our continuity result (condition {\bf (C)}) and we state our result (Theorem \ref{thm:equality}). In the rest of the section we prove our statement and provide several examples satisfying our required assumptions. 

\bigskip
To finish this introduction, let us discuss some points, which are left as future research. From \cite{cont:16}, we know that the functional It\^o calculus is also valid for more general semimartingales $X$. In particular continuity is not really relevant in this framework. However the presence of jumps requires a very careful discussion about the possibility of jumping inside the singularity set of $\xi$ (see \cite{popi:16} for the Markovian case). This is the reason why we impose the continuity of $X$. Moreover to avoid again very technical arguments, we do not consider locally smooth functionals (see \cite[Definition 5.2.10]{cont:16}).

\section{Setting and known results} \label{sect:setting}

We consider a filtered probability space $(\Omega,\cF,\Prb,\bF = (\cF_t)_{t\geq 0})$. The filtration is assumed to be complete and right continuous. Note that all martingales have right continuous modifications in this setting and we will always assume that we are taking the right continuous version, without any special mention. 
We assume that $(\Omega,\cF,\Prb,\bF = (\cF_t)_{t\geq 0})$ supports a $d$-dimensional Brownian motion $W$. 
In this paper for a given $T\geq 0$, we denote:
\begin{itemize}
\item $\cD$ (resp. $\cD(0,T)$): the set of all predictable processes on $\bR_+$ (resp. on $[0,T]$). $L^2_{loc}(W)$ is the subspace of $\cD$ such that for any $t\geq 0$ a.s.
$$\int_0^t |Z(s)|^2 ds < +\infty.$$
\item $ \cM_{loc}$: the set of c\`adl\`ag local martingales orthogonal to $W$.  
If $M \in \cM_{loc}$ then
$$[ M, W^i](t) =0, 1\leq i \leq k.$$ 
\item $\cM$ is the subspace of $ \cM_{loc}$ of martingales.
\end{itemize}
On $\bR^d$, $|.|$ denotes the Euclidean norm and $\bR^{d\times d'}$ is identified with the space of real matrices with $d$ rows and $d'$ columns. If $z \in  \bR^{d\times d'}$, we have $|z|^2 = \mbox{trace}(zz^*)$. 

Now to define the solution of our BSDE, let us introduce the following spaces for $p\geq 1$.
\begin{itemize}
\item $\bD^p(0,T)$ is the space of all adapted \cad processes $X$ such that
$$\E \left(  \sup_{t\in [0,T]} |X(t)|^p \right) < +\infty.$$
For simplicity, $X_* = \sup_{t\in [0,T]} |X(t)|$.
\item $\bH^p(0,T)$ is the subspace of all processes $X\in \cD(0,T)$ such that
$$\E \left[ \left( \int_0^T |X(t)|^2 dt\right)^{\frac{p}{2}} \ \right] < +\infty.$$
\item $\bM^p(0,T)$ is the subspace of $\cM$ of all martingales such that
$$\E \left[ \left( [ M ](T) \right)^{\frac{p}{2}}\right] < +\infty.$$
\item $\bS^p(0,T) = \bD^p(0,T) \times \bH^p(0,T)  \times \bM^p(0,T)$.
\end{itemize}
If $M$ is a $\bR^d$-valued martingale in $\cM$, the bracket process $[ M ](t)$ is
$$[ M ](t) = \sum_{i=1}^d [ M^i ](t),$$
where $M^i$ is the $i$-th component of the vector $M$. 

We consider the BSDE \eqref{eq:BSDE_jumps}
\begin{equation*} 
Y(t) = \xi + \int_t^T f(s,Y(s), Z(s)) ds  -\int_t^T Z(s) dW(s) - \int_t^T dM(s).
\end{equation*}
Here, the random variable $\xi$ is $\cF_T$-measurable with values in $\bR$ and the generator $f : \Omega \times [0,T] \times \bR \times \bR^{d} \to \bR$ is a random function, measurable with respect to $\mbox{Prog }\times \mathcal{B}(\bR)\times \mathcal{B}(\bR^{d})$ where $\mbox{Prog}$ denotes the sigma-field of progressive subsets of $\Omega \times [0,T]$. The unknowns are $(Y,Z,M)$ such that
\begin{itemize}
\item $Y$ is progressively measurable and \cad with values in $\bR$;
\item $Z \in L^2_{loc}(W)$, with values in $\bR^{d}$;
\item $M \in \cM_{loc}$ with values in $\bR$.
\end{itemize}
For notational convenience we will denote $f^0(t) = f(t,0,0)$. 

\subsubsection*{Assumptions}

\begin{itemize}
\item $\xi$ and $f^0$ are non negative and $\Prb(\xi = +\infty) > 0$. $\cS$ is the set of singularity: 
$$\cS = \{\xi=+\infty \}.$$
\item The function $y\mapsto f(t,y,z)$ is continuous and monotone: there exists $\chi \in \bR$ such that a.s. and for any $t \in [0,T]$ and $z \in \bR^k$
\begin{equation}\label{eq:f_mono} \tag{A1}
(f(t,y,z)-f(t,y',z))(y-y') \leq \chi (y-y')^2.
\end{equation}
\item For every $n> 0$ the function
\begin{equation}\label{eq:f_growth_y_t} \tag{A2}
\sup_{|y|\leq n} |f(t,y,0)-f^0_t| \in L^1((0,T) \times \Omega).
\end{equation}
\item $f$ is Lispchitz in $z$, uniformly w.r.t. all parameters: there exists $L > 0$ such that for any $(t,y)$, $z$ and $z'$: a.s.
\begin{equation}\label{eq:f_lip_z} \tag{A3}
|f(t,y,z)-f(t,y,z')| \leq L |z-z'|.
\end{equation}
\end{itemize}
Note that no assumption on $f^0$ (expect non negativity) is required. Conditions \eqref{eq:f_mono}-\eqref{eq:f_lip_z} will ensure existence and uniqueness of the solution for a version of BSDE \eqref{eq:BSDE_jumps}, where the terminal condition $\xi$ is replaced by $\xi \wedge n$ and where the generator $f$ is replaced by $f_n = f-f^0 + (f^0\wedge n)$ for some $n>0$ (see BSDE \eqref{eq:trunc_BSDE} below). We obtain the minimal supersolution (see Theorem \ref{thm:exists_super_sol}) with singular terminal condition $\xi$ by letting the truncation $n$ tend to $\infty$. To ensure that in the limit (when $n$ goes to $\infty$) the solution component $Y$ attains the value $\infty$ on $\cS$ at time $T$ but is finite before time $T$, we suppose that 
\begin{itemize}
\item There exists a constant $q > 0$ and a positive process $a$ such that for any $y \geq 0$
\begin{equation}\label{eq:f_upper_bound} \tag{A4}
f(t,y,z)\leq - (a(t)) y|y|^{q} + f(t,0,z).
\end{equation}

\end{itemize}
Moreover, in order to derive the a priori estimate, the following assumptions will hold.
\begin{itemize}
\item There exists some $\ell > 1$ such that
\begin{equation}\label{eq:alpha_gamma} \tag{A5}
\E \int_0^T \left[ \left( \frac{1}{qa(s)}\right)^{\frac{\ell}{q}}+ \left(  f^0(s)\right)^{\ell}\right] ds < +\infty.
\end{equation}
\end{itemize}
\begin{defi}
The generator $f$ satisfies Condition {\rm \textbf{(A)}} if all assumptions \eqref{eq:f_mono}--\eqref{eq:alpha_gamma} hold.
\end{defi}

\begin{ex}[Toy example] \label{ex:toy_ex}
{\rm 
The function $f(y) = -y|y|^q$ satisfies all previous conditions. It corresponds to generator \eqref{eq:example_control} with $\alpha_t = (1/q)^{1/q}$ and $\gamma_t = 0$. 

\hfill $\diamond$
}
\end{ex}

\begin{rem}
In \cite{krus:popi:16b} or in \cite{popi:16}, we consider some weaker integrability conditions on $f^0$ (see  {\rm (A6)} in \cite{krus:popi:16b} and {\rm (A$6^*$)} and {\rm (A8)} in \cite{popi:16}). These weak hypotheses can be also assumed here. But since it is not the core of this paper, we work under this stronger (but easier to check) condition {\rm (A5)} on $f^0$.
\end{rem}

\subsection{Known results}

In \cite{krus:popi:16,krus:popi:17}, we proved that if $\xi \in L^p(\Omega)$, for some $p >1$, then under Conditions {\rm {\bf (A)}} there exists a unique solution $(Y,Z,M)$ in $\bS^p(0,T)$ to the BSDE \eqref{eq:BSDE_jumps}. In \cite{krus:popi:16b}, the following result is proved. 
\begin{prop}[Theorem 1 in \cite{krus:popi:16b}]\label{thm:exists_super_sol} 
Under Condition {\rm {\bf (A)}} there exists a process $(Y,Z,M)$ such that 
\begin{itemize}
\item $(Y,Z,M)$ belongs to $\bS^\ell(0,t)$ for any $t < T$.  
\item A.s. for any $t\in [0,T]$, $Y_t \geq 0$.
\item For all $0\leq s \leq t < T$:
\end{itemize}
\begin{equation}\label{eq:bsde_dynamic}
Y(s) = Y(t) + \int_s^t f(t,Y(r),Z(r)) dr-\int_s^t Z(r) dW(r) + M(t) -M(s).
\end{equation}
\begin{itemize}
\item If the filtration $\bF$ is left-continuous at time $T$, $(Y,Z,M)$ is a super-solution in the sense that a.s. \eqref{eq:term_cond_super_sol} holds:
\begin{equation*}
\liminf_{t\to T} Y(t) \geq \xi.
\end{equation*}
\end{itemize}
Any process $(\tilde Y, \tilde Z, \tilde M)$ satisfying the previous four items is called \textbf{super-solution} of the BSDE \eqref{eq:BSDE_jumps} with singular terminal condition $\xi$.
\end{prop}
 In \cite{krus:popi:16b}, a key point in the construction of the solution is the following a priori estimate:
\begin{eqnarray}\label{eq:a_priori_estimate}
Y(t) & \leq & \frac{K_{\ell,L}}{(T-t)^{1+1/q}} \left\{\E \left( \ \int_t^{T} \left[ \left(\frac{1}{qa(s)}\right)^{1/q} + (T-s)^{1+1/q} f^0(s) \right]^{\ell} ds \bigg| \cF_t\right) \right\}^{1/\ell} \\ \nonumber
& = & \frac{K_{\ell,L}}{(T-t)^{1+1/q}} \Gamma(t)
\end{eqnarray}
where $K_{\ell,L}$ is a non negative constant depending only on $\ell$ and $L$ and this constant is a non decreasing function of $L$ and a non increasing function of $\ell$. Condition \eqref{eq:alpha_gamma} implies that a.s. $Y_t < +\infty$ on $[0,T)$. 

\begin{rem}
$\ $
\begin{itemize}
\item The constants $K_{\ell,L}$ and $\ell > 1$ come from the growth condition on $f$ w.r.t. $z$. 
\item If $f(y)=-y|y|^q$, $a(t) =1$, $L=0$, and we obtain as in \cite{popi:06}:
$$Y_t \leq \left( \frac{1}{q(T-t)} \right)^{1/q}.$$
\item The solution $(Y,Z,M)$ obtained by approximation is minimal, that is if $(\widetilde Y, \widetilde Z,\widetilde M)$ is another non negative super-solution, then for all $t \in [0,T]$, $\Prb$-a.s. $\widetilde Y_t \geq Y_t$.
\end{itemize}
\end{rem}

\vspace{0.5cm}
Now we give the main ideas of the proof of the existence result (Theorem \ref{thm:exists_super_sol}). It is important to study the behaviour of $Y$ in the next sections. The approach in \cite{krus:popi:16b} is to approximate our BSDE by considering a terminal condition of the form $\xi^n:=\xi\wedge n$ and observe asymptotic behaviour. In the rest of the paper, $(Y^n,Z^n,M^n)$ will be the solution of the truncated BSDE:
\begin{eqnarray}\label{eq:trunc_BSDE}
Y^n(t) & = & \xi\wedge n + \int_t^T f_n(s,Y^n(s),Z^n(s)) ds -\int_t^T Z^n(s) dW(s)  - \int_t^T dM^n(s).
\end{eqnarray}
Here $f_n(t,y,z)$ is the generator obtained by the truncation on $f^0$:
\begin{equation} \label{eq:truncated_gene}
f_n(t,y,z) = (f(t,y,z)-f^0(t)) + (f^0(t) \wedge n).
\end{equation}
Existence and uniqueness of $(Y^n,Z^n,M^n)$ comes from Theorem 2 in \cite{krus:popi:16}. Moreover using comparison argument (see \cite{krus:popi:16} or \cite{quen:sule:13}) we can obtain for $m \leq n$: $0\leq Y^m(t) \leq Y^{n}(t) $. And for any $n$, $Y^n$ satisfies Estimate \eqref{eq:a_priori_estimate} (Proposition 2 in \cite{krus:popi:16b}). 
This allows us to define $Y$ as the limit of the increasing sequence
$(Y^n_t)_{n\geq 1}$: $$\forall\, t\in[0,T],\quad Y(t) :=\lim_{n\rightarrow \infty} Y^n(t).$$
Proposition 3 in \cite{krus:popi:16b} shows that there exists a constant $C$ such that for any $0 < t < T$
\begin{eqnarray}\label{eq:approx_estim_L2_main_ineq}
&&\E\left[\sup_{0\leq s\leq t} |Y^n(s)-Y^m(s)|^\ell+\left( \int_0^t |Z^n(s)-Z^m(s)|^2 ds \right)^{\ell/2} + [M^n - M^m](t)^{\ell/2}\right]\\ \nonumber
&& \\ \nonumber
&& \qquad \qquad \leq C\E\left[ |Y^n(t)-Y^m(t)|^\ell\right] + C \E \int_0^t |f^0(s)\wedge n - f^0(s)\wedge m|^\ell ds. 
\end{eqnarray}
Since $Y^n(t)$ converges to $Y(t)$ almost surely, with the a priori estimate \eqref{eq:a_priori_estimate}, Condition \eqref{eq:alpha_gamma} and Inequality \eqref{eq:approx_estim_L2_main_ineq}, thanks to the dominated convergence theorem, we can deduce that for every $\varepsilon>0$, $(Y^n,Z^n,M^n)_{n\geq1}$ converges to $(Y,Z,M)$ in $\bS^\ell(0,T-\varepsilon)$. The limit $(Y,Z,M)$ satisfies for every $0\leq t<T$, for all $0\leq s \leq t$ the dynamics \eqref{eq:bsde_dynamic} and $Y$ satisfies Inequality \eqref{eq:a_priori_estimate}. Note that all these results are obtained without the left-continuity assumption on the filtration $\bF$. But since the solution $(Y,Z,M)$ satisfies the dynamic \eqref{eq:bsde_dynamic} only on $[0,T-\eps]$ for any $\eps > 0$, we cannot derive directly the existence of a left limit at time $T$ for $Y$. 

Note that on $Z$ we also have a stronger integrability result:
\begin{prop} \label{prop:sharp_estim_Z_U}
Under Condition {\bf (A)}, there exists a constant $C$ independent of $n$ such that the process $Z^n$ satisfies:
$$\E\left[ \int_0^T (T-s)^{\rho} \left( |Z^n(s)|^2  \right) ds  \right]^{\ell/2} \leq C.$$
The constant $\rho$ satisfies $\displaystyle \rho > \frac{2}{q} + 2\left( 1 -\frac{1}{\ell} \right) = \frac{2}{q} + \frac{2}{\ell^*}$, if $\ell^*$ is the H\"older conjugate of $\ell$.
\end{prop}
In the sequel we assume that $q > 2$. W.l.o.g. we can assume that $1< \ell < \frac{2q}{2+q}$ in condition \eqref{eq:alpha_gamma} such that $\frac{2}{q} + 2\left( 1 -\frac{1}{\ell} \right)  < 1$ and we can choose $\rho < 1$ in the above proposition.
In particular if the generator is $f(y) = - y|y|^q$ (Example \ref{ex:toy_ex}), then we can take  $q > 2$ and $\ell =1$, which was supposed in \cite{popi:06}. In this particular case, the constant $C$ is explicitely given by: $C = 16\left(1/q\right)^{2/q}$.

\subsection{Functional It\^o calculus}

We adopt the notations and the setting developed in \cite{cont:four:13,cont:16}. We simply copy the main definitions; all details can be found in these two works. 

Consider $D([0,T],\bR^d\times S^+_d)$, the space of all c\`adl\`ag functions defined on $[0,T]$ with values in $\bR^d\times S^+_d$, where $S^+_d$ is the set of positive $d\times d$ matrices.
For a path $\omega \in D([0,T],\bR^d)$, $\omega_t=\omega (t \wedge \cdot)$ is the path stopped at time $t$ and $\omega_{t-} = \omega \mathbf{1}_{[0,t[} + \omega(t-) \mathbf{1}_{[t,T]}$, where $\omega(t-)$ is the left limit of $\omega$ at time $t$. 

The space of stopped paths $\Upsilon$ is defined as the quotient space of $[0,T] \times D([0,T],\bR^d\times S_d^+)$ by the equivalence relation:
$$(t,\omega) \sim (t',\omega') \Leftrightarrow (t=t'\ \mbox{and} \ \omega_t = \omega_{t'}).$$
This space is endowed with the distance 
$$d_\infty((t,\omega) ,(t',\omega')) = \|\omega_t - \omega_{t'}\|_\infty + |t-t'|.$$
And $(\Upsilon,d_\infty)$ is a metric space and a closed subset of $([0,T] \times D([0,T],\bR^d\times S^+_d); \|\cdot\|_\infty)$.
The notion of non-anticipative, continuous, left-continuous and boundedness-preserving functionals is defined in \cite[Definitions 2.1, 2.3, 2.4 and 2.5]{cont:four:13} and we denote by $\bC^{0,0}([0,T))$ (resp. $\bC^{0,0}([0,T))$, resp. $\bB([0,T))$) the set of continuous (resp. left-continuous, resp. boundedness-preserving) functions on $\Upsilon$. In the sequel we assume that $F :\Upsilon \to \bR$ is a non-anticipative functional with {\it predictable dependence} with respect to $v$:
\begin{equation} \label{eq:predictable_dep_cond}
\forall (t,x,v) \in \Upsilon, \quad F(t,x,v)=F(t,x_t,v_{t-}).
\end{equation}


Let us briefly recall the definition of the horizontal and vertical derivatives. Let $F :\Upsilon \to \bR$ be a non-anticipative functional. The horizontal derivative $\cD F$ of $F$ at $(x,v) \in D([0,t],\bR^d\times S^+_d)$ is the limit (if it exists) 
$$\cD F(t,x,v) = \lim_{h \downarrow 0} \frac{F(t+h,x_t,v_t)-F(t,x_t,v_t)}{h}.$$ 
If the limit exists for all $(x,v)$, the map $\cD : \Upsilon\to \bR$ defines a non-anticipative functional $\cD F$, called the {\it horizontal derivative} (see \cite[Definition 3.1]{cont:four:13}).

Now from \cite[Definition 3.2]{cont:four:13}, $F$ is {\it vertically differentiable} at $(t,x,v)$ if the map defined on $\bR^d$ by $e\mapsto F(t,x_t + e\mathbf{1}_{[t,T]},v_t)$ is differentiable at 0. The vertical derivative $\nabla_\omega F$ at $(t,x,v)$ is the gradient of the previous map:
$$\nabla_\omega F(t,x,v) = \left( \lim_{h\to 0} \frac{F(t,x_t + h e_i\mathbf{1}_{[t,T]},v_t)-F(t,x_t,v_t)}{h}, \ i=1,\ldots,d\right).$$
Let us recall \cite[Definition 3.6]{cont:four:13}:
\begin{defi}[$\bC^{1,k}$ functionals]
Define $\bC^{1,k}([0,T))$ as the set of left-continuous functionals $F \in \bC^{0,0}_l$ such that 
\begin{itemize}
\item $F$ admits a horizontal derivative $\cD F(t,\omega)$ for all $(t,\omega) \in \Lambda_T$, and the map $\cD F(t,\cdot) : (D([0,T],\bR^d),\|.\|_\infty) \to \bR$ is continuous for each $t \in [0,T[$;
\item $F$ is $k$ times vertically differentiable with $\nabla^j_\omega F \in \bC^{0,0}_l$.
\end{itemize}
We define $\bC^{1,k}_b([0,T))$ as the set of functionals $F \in \bC^{1,k}$ such that $\cD F$, $\nabla_\omega F,\ldots,\nabla^k_\omega F$ belong to $\bB(\Lambda_T)$.
\end{defi}

In the sequel we will use the change of variable formulas.

\noindent {\bf Theorem 4.1 in \cite{cont:four:13}.}
{\it Let $F \in \bC^{1,2}_b$ verifying \eqref{eq:predictable_dep_cond}. Let $X$ be a continuous $\bR^d$-valued semimartingale with absolutely continuous quadratic variation 
$$[X](t) = \int_0^t A(u) du,$$
where $A$ is an $S^+_d$-valued process. Then for $t \in [0,T]$
\begin{eqnarray} \label{eq:functional_Ito_formula}
&& F(t,X_t,A_t) = F(0,X_0,A_0)  + \int_0^t \cD F(u,X_u,A_u)du \\ \nonumber
&& \qquad + \int_0^t \nabla_\omega F(u,X_u,A_u)dX(u) +  \frac{1}{2} \int_0^t \tr \left( \nabla^2_\omega F(u,X_u,A_u)d[X](u) \right) .
\end{eqnarray}
}
This result implies in particular that $\cX = F(\cdot,X,A)$ is a continuous semimartingale for any $F \in \bC^{1,2}_b$. In the sequel, we need some integrability properties of $\cX$. Let us recall that the classical norm on semimartingales is defined in \cite{dell:meye:80}, Section VII.3 (98.1)-(98.2) or \cite{prot:04}, Section V.2. Nevertheless this norm is not sufficient in our case and we follow the ideas of \cite[Section 7.5]{cont:16}. For $p \geq 1$, we define $\cA^p(\bF)$ as the set of continuous $\bF$-predictable absolutely continuous processes $H=H(0) + \int_0^\cdot h(t) dt$ with finite variation such that 
$$\|H\|^p_{\cA^p} = \bE\left(  |H(0)|^p + \int_0^T |h(t)|^p dt \right) < +\infty.$$
We consider the direct sum
$$\cS^{p} = \bM^p(0,T)\oplus \cA^p(\bF).$$
Any process $S \in \cS^{p}$ is an $\bF$-adapted special semimartingale with a unique decomposition $S=M+H$, where $M \in  \bM^p(0,T)$ with $M(0)=0$ and $H \in \cA^p(\bF)$ with $H(0)=0$. Let us remark that by Jensen's inequality, the norm defined on $\cS^{p}$ is stronger than the norm of semimartingales defined in \cite{dell:meye:80}. Moreover if $S \in \cS^{p}$, then $S \in \bD^p(0,T)$ by the Burkh\"older-Davis-Gundy inequality. The interested reader can find in \cite[Chapter 7]{cont:16} how the vertical and horizontal derivatives can be defined on this space $\cS^{p}$.

%

\section{Continuity at time $T$ in the non-Markovian setting} \label{sect:continuity_Y_Markov}

First note that we do not impose any further condition on the generator. In the sequel we assume that $X$ is the solution of the SDE \eqref{eq:SDE}
\begin{equation} \label{eq:SDE}
X(t)=\zeta(t)+\int_0^t b(s,X_s)\dd s +\int_0^t \sigma(s,X_s)\dd W(s)
\end{equation}
The coefficients $b(\cdot,\cdot,\phi):\Omega \times [0,T] \to \bR^d$ and $\sigma(\cdot,\cdot,\phi):\Omega \times [0,T] \to \bR^{d\times d}$ are defined for every continuous function $\phi$ and satisfy the standard conditions:
\begin{itemize}
\item $b$, $\sigma$ are Lipschitz continuous w.r.t. $\phi$ uniformly in $t$ and $\omega$, i.e. there exists a constant $K_{b,\sigma}$ such that for any $(\omega,t) \in \Omega \times [0,T]$, for any $\phi$ and $\psi$ in $C([0,T];\bR^d)$: a.s.
\begin{equation*}
|b(t,\phi)-b(t,\psi)| + |\sigma(t,\phi)-\sigma(t,\psi)|  \leq K_{b,\sigma} \|\phi_t-\psi_t\|_{\infty}.
\end{equation*}
\item $b$ and $\sigma$ growth at most linearly:
\begin{equation*}
|b(t,0)| + |\sigma(t,0)|  \leq C_{b,\sigma}.
\end{equation*}
\item $\zeta$ is a progressively measurable continuous stochastic process such that for some $\varrho \geq 0$, $\zeta \in \bD^\varrho(0,T)$.  
\end{itemize}
Recall that $\phi_t$ is the stopped path of $\phi$, which implies that 
$$\|\phi_t-\psi_t\|_{\infty} = \sup \{|\phi(u)-\psi(u)|, \ 0 \leq u \leq t\}.$$
Let us emphasize that $X$ is not a Markovian process since the drift and the volatility matrix may depend on the whole trajectory of $X$. 
Under those assumptions, the forward SDE \eqref{eq:SDE} has a unique strong continuous solution $X$ (see \cite[Theorem 3.17]{pard:rasc:14}), such that 
\begin{equation} \label{eq:Dp_estim_X}
\bE \left[ \sup_{t\in [0,T]} |X(t)|^\varrho \right] \leq C_{\varrho}.
\end{equation} 
To lighten the notation, the dimensions of $X$ and of the Brownian motion are the same. But this condition is not crucial and we can also work with different dimensions. The process $X$ is a continuous semimartingale with
$$[X](t) = \int_0^t \sigma(s,X_s)\sigma^*(s,X_s) ds = \int_0^t A(s) ds.$$
We assume that {\bf Condition (C)} holds, namely:
\begin{itemize}
\item[C1.] There exists a measurable function $\Phi : \bR \to [0,+\infty]$ and $F \in \bC^{1,2}_b$ such that 
$$\xi = \Phi(F(T,X_T,A_T)).$$
We denote $\cR=\{ \Phi < +\infty\}$ which is supposed to be an open subset of $\bR$ and we suppose that $\bP(\xi=\infty) > 0$. 
\item[C2.] For any compact set $\cK \subset \cR$, $\E (\xi \ind_{\cK}(F(T,X_T,A_T))) < +\infty$.
\item[C3.] $F(\cdot,X,A)$ is in $\cS^{p}$ for $p =  \frac{q+1}{q}\ell^*$, where $\ell^*$ is the H\"older conjugate of the constant $\ell > 1$ of Condition \eqref{eq:alpha_gamma}.
\item[C4.] $\nabla_\omega F$ is in $\bD^{\ell^*}(0,T)$ and $\varrho$ in \eqref{eq:Dp_estim_X} is equal to $\ell^*$.
\end{itemize}
Recall that (C3) implies that $F(\cdot,X,A)$ is in $\bD^p(0,T)$
$$\E \left[\sup_{t\in [0,T]} |F(t,X_t,A_t)|^p \right] < +\infty,$$
and that $\sigma(\cdot,X) \nabla_\omega F(\cdot,X,A)$ belongs to $\bH^p(0,T)$
\begin{equation} \label{eq:Lp_estim_vert_deriv}
\E \left[ \left(  \int_0^T\left[(\nabla_\omega F(s,X_s,A_s))^* A(s) \nabla_\omega F(s,X_s,A_s) \right] \dd s \right)^{p/2}\right]< +\infty.
\end{equation}
Moreover the It\^o formula \eqref{eq:functional_Ito_formula}, together with \eqref{eq:SDE}, implies that 
$$ F(t,X_t,A_t) = F(0,X_0,A_0)  + \int_0^t \Theta_1(u)du  + \int_0^t \nabla_\omega F(u,X_u,A_u)\sigma(u,X_u)dW(u) ,$$
with
\begin{eqnarray*}
 \Theta_1(s)& =&\bigg\{  \cD F(s,X_s,A_s) + \nabla_\omega F(s,X_s,A_s) b(s,X_s)  \\ \nonumber
&& \qquad \qquad\qquad \qquad \qquad \left. +  \frac{1}{2} \tr \left( \nabla^2_\omega F(s,X_s,A_s) A(s) \right) \right\} .
\end{eqnarray*}
From (C3), we obtain that 
$$\bE \int_0^T |\Theta_1(s)|^p ds < +\infty.$$
\begin{rem}
Recall that if $q > 2$, in order to apply Proposition \ref{prop:sharp_estim_Z_U}, we have chosen $\ell < \dfrac{2q}{q+2}$. Hence 
$$p = \frac{q+1}{q} \ell^* > 2 \frac{q+1}{q-2} = 2 +\frac{6}{q-2}.$$
In particular if $q$ is close to 2, $p$ is large.

Instead of condition (C4), we may assume that $\nabla_\omega F$ is also in $\cS^p$, following the idea of \cite[Section 7.5]{cont:16}.
\end{rem}

Let us state our main result.
\begin{thm} \label{thm:equality}
Under the hypotheses {\rm \textbf{(A)}} and {\rm {\bf (C)}}, with $q > 2$, the minimal supersolution $Y$ satisfies a.s. 
$$\liminf_{t\to T} Y_t = \xi.$$
\end{thm}

\subsection{Proof of Theorem \ref{thm:equality}} \label{sect:equality}

Let $(Y^n,Z^n,M^n)$ be the solution of the BSDE \eqref{eq:trunc_BSDE} with terminal condition $\xi \wedge n$ and generator $f_n$
$$Y^n(t) = \xi\wedge n + \int_t^T f_n(s,Y^n(s),Z^n(s)) \dd s -\int_t^T Z^n(s) \dd W(s) - M^n(T) + M^n(t).$$
Let us emphasize that the process $Y^n$ is bounded on $\Omega \times [0,T]$. 

Let $\phi:\bR \to \bR_+$ be a $C^\infty$-function with compact support included in $\cR$.
We apply It\^o's formula to the process $Y^n\phi(F(\cdot,X,A))$ between $0$ and $t$: 
\begin{eqnarray*}
&&Y^n(t) \phi(F(t,X_t,A_t)) = Y^n(0) \phi(F(0,X_0,A_0))+\int_0^t \phi(F(s,X_s,A_s))Z^n(s) dW(s) \\
&&\quad + \int_0^t Y^n(s)  \phi^\prime (F(s,X_s,A_s)) \nabla_\omega F(s,X_s,A_s) \sigma(s,X_s)dW(s)\\
&&\quad + \int_0^t  \phi (F(s,X_s,A_s)) dM^n(s)\\
&&\quad + \int_0^t  \phi^\prime (F(s,X_s,A_s)) \nabla_\omega F(s,X_s,A_s) \sigma(s,X_s)Z^n(s) ds  \\
&&\quad - \int_0^t f_n(s,Y^n(s),Z^n(s)) \phi(F(s,X_s,A_s)) ds\\
&&\quad + \int_0^t Y^n(s) \phi^\prime(F(s,X_s,A_s)) \Theta_1(s) ds +\frac{1}{2} \int_0^tY^n(s) \phi^{\prime\prime}(F(s,X_s,A_s)) \Theta_2(s) ds.
\end{eqnarray*}
with
\begin{eqnarray*}
\Theta_2(s) &=&  \frac{1}{2}  \left( \nabla_\omega F(s,X_s,A_s) \sigma(s,X_s) \right)  \left( \nabla_\omega F(s,X_s,A_s) \sigma(s,X_s) \right)^* .
\end{eqnarray*}

Now we decompose the quantity with the generator $f_n$ as follows:
\begin{eqnarray*} 
&&\int_0^t \phi(F(s,X_s,A_s)) f_n(s,Y^n(s),Z^n(s)) ds \\ \nonumber
&& \quad = \int_0^t \phi(F(s,X_s,A_s)) (f(s,Y^n(s),0) -f^0(s)) ds  \\ \nonumber
&& \qquad + \int_0^t \phi(F(s,X_s,A_s)) (f^0(s)\wedge n) ds \\ \nonumber
&&\qquad + \int_0^t \phi(F(s,X_s,A_s)) \left( f(s,Y^n(s),Z^n(s)) - f(s,Y^n(s),0)\right) ds \\ \nonumber
&& \quad =\int_0^t \phi(F(s,X_s,A_s)) (f(s,Y^n(s),0) -f^0(s)) ds \\ \nonumber
&& \qquad +\int_0^t \phi(F(s,X_s,A_s)) (f^0(s)\wedge n) ds \\ \nonumber
&&\qquad + \int_0^t \phi(F(s,X_s,A_s)) \zeta^n(s) Z^n(s) ds
\end{eqnarray*}
where $\zeta^n_s$ is a $d$-dimensional random vector defined by: for $i=1,\ldots,d$
$$\zeta^{i,n}(s) = \frac{\left( f(s,Y^n(s),Z^n(s)) - f(s,Y^n(s),0)\right) }{Z^{i,n}(s)} \ind_{Z^{i,n}(s)\neq 0}.$$
From Condition \eqref{eq:f_lip_z}, $|\zeta^n(s)|\leq K$. Hence we obtain
\begin{eqnarray} \label{eq:applied_func_Ito_formula}
&&Y^n(t) \phi(F(t,X_t,A_t)) = Y^n(0) \phi(F(0,X_0,A_0))\\ \nonumber
&& \quad +\int_0^t \phi(F(s,X_s,A_s)) \left[ Z^n(s) dW(s) + dM^n(s)\right]\\  \nonumber
&&\quad + \int_0^t Y^n(s)  \phi^\prime (F(s,X_s,A_s)) \nabla_\omega F(s,X_s,A_s) \sigma(s,X_s)dW(s)\\ \nonumber
&&\quad + \int_0^t  \Psi^n(s)Z^n(s) ds - \int_0^t \phi(F(s,X_s,A_s)) (f^0(s)\wedge n) ds \\ \nonumber
&&\quad - \int_0^t \phi(F(s,X_s,A_s)) (f(s,Y^n(s),0) -f^0(s)) ds  \\ \nonumber
&&\quad + \int_0^t Y^n(s) \left[ \phi^\prime(F(s,X_s,A_s)) \Theta_1(s) +  \phi^{\prime\prime}(F(s,X_s,A_s)) \Theta_2(s)\right] ds 
\end{eqnarray}
with
\begin{eqnarray*}
\Psi^n(s) &=&  \phi^\prime (F(s,X_s,A_s)) \nabla_\omega F(s,X_s,A_s) \sigma(s,X(s)) +\phi(F(s,X_s,A_s)) \zeta^n(s) .
\end{eqnarray*}
Recall that for a fixed $n$, $Y^n$ is bounded, $(Z^n,M^n)$ belong to $\bH^p(0,T) \times \bM^p(0,T)$ for any $p \geq 1$. From Condition (C3) on $F(\cdot,X,A)$, taking the expectation in \eqref{eq:applied_func_Ito_formula} leads to: for $t \in [0,T]$
\begin{eqnarray} \label{eq:expect_func_Ito_formula}
&&\bE \left[Y^n(T) \phi(F(T,X_T,A_T)) \right]= \bE \left[Y^n(t) \phi(F(t,X_t,A_t)) \right] \\ \nonumber
&&\quad - \bE \left[\int_t^T \phi(F(s,X_s,A_s)) (f^0(s)\wedge n) ds \right] \\ \nonumber
&&\quad - \bE \left[\int_t^T \phi(F(s,X_s,A_s)) (f(s,Y^n(s),0) -f^0(s)) ds \right] \\ \nonumber
&&\quad +\bE \left[ \int_t^T  Y^n(s) \left[ \phi^\prime(F(s,X_s,A_s)) \Theta_1(s) +  \phi^{\prime\prime}(F(s,X_s,A_s)) \Theta_2(s)\right] ds \right]\\ \nonumber
&&\quad +\bE \left[ \int_t^T  \Psi^n(s)Z^n(s) ds \right] .
\end{eqnarray}

From the assumptions (C1) and (C2) on $\xi=\Phi(F(T,X_T,A_T))$, we have for any $n$:
\begin{equation} \label{eq:control_cont_T}
\E (Y^n(T) \phi(F(T,X_T,A_T)) \leq \E [ \Phi(F(T,X_T,A_T))\phi(F(T,X_T,A_T))] < +\infty.
\end{equation}
From the a priori estimate \eqref{eq:a_priori_estimate}, Assumption \eqref{eq:alpha_gamma} and from the boundedness of $\phi$, for any $t < T$
\begin{equation} \label{eq:control_cont_t}
\E (Y^n(t) \phi(F(t,X_t,A_t))) \leq \frac{1}{(T-t)^{1/q+1}} \E (\Gamma(t)\phi(F(t,X_t,A_t))) < +\infty
\end{equation}
Since $\phi$ is bounded and $f^0 \in L^1((0,T)\times \Omega)$ (Condition \eqref{eq:alpha_gamma}):
\begin{equation} \label{eq:control_cont_0}
\E \int_0^T \phi(F(s,X_s,A_s)) (f^0_s\wedge n) ds \leq C.
\end{equation}
%

Now we treat the two terms in \eqref{eq:expect_func_Ito_formula} containing $Y^n$. Firstly, by condition \eqref{eq:f_upper_bound} remark that:
\begin{eqnarray} \label{eq:control_cont_1}
&&  - \int_0^t \phi(F(s,X_{s},A_s)) \left( f(s,Y^n(s),0)-f(s,0,0) \right) ds\\ \nonumber
&& \qquad  \geq \int_0^t \phi(F(s,X_{s},A_s)) a(s) |Y^n(s)|^{1+q} ds.
\end{eqnarray}
Secondly with H\"older's inequality we obtain for $j=1$ or 2:
\begin{eqnarray*}
&& \int_0^t |Y^n(s) \phi^{(j)}(F(s,X_{s},A_s))\Theta_j(s)| ds  \leq  \left[  \int_0^t  \phi(F(s,X_{s},A_s)) a(s) |Y^n(s)|^{1+q}  ds \right]^{\frac{1}{q+1}}\\ 
&& \hspace{2cm} \times \left[ \int_0^t a(s)^{-1/q} \phi(F(s,X_{s},A_s))^{-1/q}|\phi^{(j)}(F(s,X_{s},A_s))|^{\frac{q+1}{q}} |\Theta_j(s)|^{\frac{q+1}{q}}  ds \right]^{\frac{q}{q+1}}.
\end{eqnarray*}
To control the second quantity, we will be more specific about the test-function $\phi$. We will assume that $\phi = \psi^\gamma$ where $\psi$ belongs to $C^\infty_b(\bR^d)$ with support in $\cR$ and $\gamma> 2(q+1)/q$. Under this setting, there exists a constant $C$ depending only on $\psi$ and $\gamma$ such that 
$$|\phi^\prime| + |\phi^{\prime\prime}| \leq C \psi^{\gamma-2}.$$ 
Thus for $\gamma > 2(q+1)/q$ and $j=1$ or 2
\begin{eqnarray*}
\phi(F(s,X_{s},A_s))^{-1/q} |\phi^{(j)} (F(s,X_{s},A_s))|^{(q+1)/q} & \leq & C \psi(F(s,X_{s},A_s))^{\gamma-2(q+1)/q},
\end{eqnarray*}
which is bounded. By condition \eqref{eq:alpha_gamma}, $a^{-1/q}$ is in $\bL^\ell(\Omega \times [0,T])$. By the assumption (C3), the quantity $|\Theta_j|^{\frac{q+1}{q}}$ is in $\bL^{\ell^*}(\Omega \times [0,T])$. We deduce that there exists a constant $C$ such that for any $t \in [0,T]$
\begin{eqnarray} \label{eq:control_cont_2}
&&  \E \left[ \int_0^t|Y^n(s) \phi^{(j)}(F(s,X_{s},A_s))\Theta_j(s)| ds  \right] \\ \nonumber
&& \qquad\leq C \left[ \E \int_0^t a(s) \phi(F(s,X_{s},A_s)) (Y^n(s))^{q+1}  ds \right]^{\frac{1}{q+1}}.
\end{eqnarray}

We use H\"older's and Young's inequalities to obtain:
\begin{eqnarray*}
&& \int_0^t |\left[ \phi^\prime (F(s,X_s,A_s)) \nabla_\omega F(s,X_s,A_s) \sigma(s,X_s) +\phi(F(s,X_s,A_s)) \zeta^n(s) \right] Z^n(s)| ds\\ \nonumber
&& \quad \leq \int_0^t \left| \Psi^n(s) \right| Z^n(s) ds \\ \nonumber
&& \quad  \leq \left[\int_0^t (T-s)^{\rho} |Z^n(s)|^2 ds\right]^{1/2}  \left[ \int_0^t \frac{|\Psi^n(s)|^2}{(T-s)^{\rho}} ds \right]^{1/2}  \\  \nonumber
&&\quad \leq \frac{1}{\ell} \left[\int_0^t (T-s)^{\rho} |Z^n(s)|^2 ds\right]^{\frac{\ell}{2}}  + \frac{1}{\ell^*} \left[ \int_0^t \frac{|\Psi^n(s)|^2}{(T-s)^{\rho}} ds \right]^{\frac{\ell^*}{2}} \\
&& \quad \leq \frac{1}{\ell} \left[\int_0^t (T-s)^{\rho} |Z^n(s)|^2 ds\right]^{\frac{\ell}{2}}  + \frac{1}{\ell^*} \left( \frac{ T^{1-\rho}}{1-\rho}  \right)^{\frac{\ell^*}{2}} \sup_{t\in [0,T]} |\Psi^n(t)|^{\ell^*}
\end{eqnarray*}
Taking the expectation and thanks to Proposition \ref{prop:sharp_estim_Z_U}, the first term on the right-hand side is bounded. For the second term, $\phi$, $\phi^\prime$ and $\zeta^n$ are bounded. From condition (C4), we deduce that there exists a constant $C$ such that for any $n$
\begin{equation} \label{eq:control_cont_3}
\E \int_0^T |\Psi^n_s Z^n_s| ds\leq C.
\end{equation}

Coming back to \eqref{eq:expect_func_Ito_formula} with $t=0$ and using \eqref{eq:control_cont_T}, \eqref{eq:control_cont_t} with $t=0$, \eqref{eq:control_cont_0}, \eqref{eq:control_cont_1}, \eqref{eq:control_cont_2} and \eqref{eq:control_cont_3}, we deduce that for any function $\phi=\psi^\gamma$ with $\gamma > 2(q+1)/q$, there exists a constant $C$ independent of $n$ such that 
\begin{equation} \label{eq:control_cont_main_term}
0\leq \E \int_0^T a(s) \phi(F(s,X_s,A_s))|Y^n(s)|^{1+q} ds \leq C < +\infty.
\end{equation}
Moreover by the monotone convergence theorem, we can pass to the limit when $n$ goes to $+\infty$ in the first four terms of \eqref{eq:expect_func_Ito_formula}. For the last two terms (containing $Y^n$ and $Z^n$), let us summarize the arguments (see details in \cite{popi:06}). Estimate \eqref{eq:control_cont_main_term} shows that the sequence $a^{\frac{1}{1+q}} \phi(F(\cdot,X,A))^{\frac{1}{1+q}} Y^n$ is bounded in $\bL^{1+q}(\Omega \times (0,T))$. Using a weak convergence result and extracting a subsequence if
necessary, and arguing as in the proof of Estimate \eqref{eq:control_cont_2}, we can pass to the limit in the term
$$\bE \left[ \int_t^T  Y^n(s) \left[ \phi^\prime(F(s,X_s,A_s)) \Theta_1(s) +  \phi^{\prime\prime}(F(s,X_s,A_s)) \Theta_2(s)\right] ds \right].$$
From Proposition \ref{prop:sharp_estim_Z_U}, there exists a subsequence, which we still denote as $(T -\cdot )^{1/(2\rho)}Z^n $, and which converges weakly in the space $\bL^2(\Omega \times (0,T))$ to a limit, and the limit is $(T -\cdot )^{1/(2\rho)}Z$, because we already know that $Z^n$ converges to $Z$ in $\bH^\ell(\Omega \times (0,T-\delta))$ for all $\delta > 0$. Let us define on $[0,T)$ the $d$-dimensional random vector $\zeta$ by: for $i=1,\ldots,d$
$$\zeta^{i}(s) = \frac{\left( f(s,Y(s),Z(s)) - f(s,Y(s),0)\right) }{Z^{i}(s)} \ind_{Z^{i}(s)\neq 0}$$
and $\Psi$ is defined as $\Psi^n$, replacing $\zeta^n$ by $\zeta$. Again $|\zeta(s)|\leq K$ and we have shown that $\Psi^n/(T-\cdot)^{1/(2\rho)}$ and $\Psi/(T-\cdot)^{1/(2\rho)}$ are in $\bL^{\ell^*}(\Omega; \bL^2(0,T))$. For any $\eps > 0$, we deduce that there exists $\delta > 0$ such that 
$$\bE \int_{T-\delta}^T \left( |\Psi^n_s Z^n_s| + |\Psi_s Z_s| \right) \dd s \leq \eps/2.$$
On the interval $[0,T-\delta]$, the sequence $(Y^n,Z^n)$ converges to $(Y,Z)$ in $\bD^\ell(0,T-\delta) \times \bH^\ell(0,T-\delta)$. Hence 
$$\lim_{n \to +\infty} \bE \int_0^{T-\delta} |\Psi^n_s Z^n_s- \Psi_s Z_s| \dd s =0.$$
In other words the sequence $\Psi^n Z^n$ converges in $\bL^1(\Omega\times (0,T))$ to $\Psi Z$ and 
$$\bE \int_0^T |\Psi(s)Z(s)| ds \leq C.$$
Passing to the limit in \eqref{eq:expect_func_Ito_formula} implies:
\begin{eqnarray} \label{eq:limit_expect_func_Ito_formula}
&&\bE \left[Y(T) \phi(F(T,X_T,A_T)) \right]= \bE \left[Y(t) \phi(F(t,X_t,A_t)) \right] \\ \nonumber
&&\quad - \bE \left[\int_t^T \phi(F(s,X_s,A_s)) f^0(s) ds \right] \\ \nonumber
&&\quad - \bE \left[\int_t^T \phi(F(s,X_s,A_s)) (f(s,Y(s),0) -f^0(s)) ds \right] \\ \nonumber
&&\quad +\bE \left[ \int_t^T  Y(s) \left[ \phi^\prime(F(s,X_s,A_s)) \Theta_1(s) +  \phi^{\prime\prime}(F(s,X_s,A_s)) \Theta_2(s)\right] ds \right]\\ \nonumber
&&\quad +\bE \left[ \int_t^T  \Psi(s)Z(s) ds \right] .
\end{eqnarray}
Estimate \eqref{eq:control_cont_main_term} also holds with $Y$, and once again from \eqref{eq:control_cont_0}, \eqref{eq:control_cont_1}, \eqref{eq:control_cont_2} and \eqref{eq:control_cont_3}, we can let $t$ go to $T$ in \eqref{eq:limit_expect_func_Ito_formula} in order to have:
$$\bE\left[ (\liminf_{t \to T} Y_t) \phi(F(T,X_T,A_T))\right] \leq  \lim_{t\to T} \bE[Y_t\phi(F(t,X_t,A_t))]= \bE[\xi \phi(F(T,X_T,A_T))] .$$
Recall that we already know that $\liminf_{t \to T} Y_t \geq \xi$ a.s. This last inequality shows that in fact a.s.
$$\liminf_{t \to T} Y_t = \xi.$$
This achieves the proof of Theorem \ref{thm:equality}.

The proof of Theorem \ref{thm:equality} shows that the limit of $Y_t$ exists in mean in the following sense: for smooth function $\phi$
$$\lim_{t \to T} \E (Y_t \phi(F(t,X_t,A_t))) =  \left\{ \begin{array}{ll}
\E (\xi \phi(F(T,X_T,A_T))) & \mbox{if  } \supp(\phi) \cap \cS = \emptyset, \\
+ \infty & \mbox{if  } \E ( \phi(F(T,X_T,A_T)) \ind_\cS) > 0.
\end{array} \right.$$

\subsection{Some examples} \label{ssect:examples}

Several examples of smooth functionals are given in  \cite{cont:four:13} or \cite{cont:16}. There are also interesting counterexamples (see \cite[Section 3.2]{cont:four:13}). 

First we can recover the Markovian case if for some smooth function $h \in C^{1,2}([0,T] \times \bR^d)$ 
$$F(t,X_t,A_t) = h(t,X(t)),$$
and if $X$ satisfies the SDE
\begin{equation} \label{eq:Markov_SDE}
X(t)=x+\int_0^t b(s,X(s))\dd s +\int_0^t \sigma(s,X(s))\dd W(s)
\end{equation}
Here the coefficients $b : [0,T] \times \bR^d \to \bR^d$ and $\sigma : [0,T] \times \bR^d\to \bR^{d\times d}$ are Lipschitz continuous w.r.t. $x$ uniformly in $t$ and $b$ and $\sigma$ growth at most linearly. Under this setting, the SDE has a unique strong continuous solution $X$ such that \eqref{eq:Dp_estim_X} holds for any $\varrho \geq 1$. Then  $\cD F(s,X_s,A_s) = \partial_t h(t,X(t))$, $\nabla_\omega F(s,X_s,A_s)= \nabla_x h(t,X(t))$ and $\nabla^2_\omega F(s,X_s,A_s)= D^2_x h(t,X(t))$, where $D^2_x$ is the Hessian matrix w.r.t. $x$. In this case Equation \eqref{eq:expect_func_Ito_formula} is the same as the classical It\^o formula used in \cite{popi:16}. If we assume that $h$ and its derivatives are of linear growth w.r.t. $x$, uniformly in time and $\omega$, then using \eqref{eq:Dp_estim_X}, the assumptions (C3) and (C4) are satisfied. 

\vspace{0.5cm}
As a second example, we consider the case where $X$ is the solution of \eqref{eq:SDE} and
$$F(t,X_t,A_t) = \int_0^t h(s,X(s)) A(s) ds$$
where $h$ is a continuous function on $[0,T] \times \bR^d$. Then $ \cD F(s,X_s,A_s) = h(s,X(s)) A(s)$, $\nabla_\omega F(s,X_s,A_s)= 0$ and Conditions (C3) and (C4) are satisfied trivially verified if $\varrho$ is sufficient large and if $h$ is of linear growth w.r.t. $x$. Moreover Equation \eqref{eq:applied_func_Ito_formula} can be simplified:
\begin{eqnarray*}
&&Y^n_t \phi(F(t,X_t,A_t)) = Y^n_0 \phi(F(0,X_0,A_0))+\int_0^t \phi(F(s,X_s,A_s)) \left[ Z^n(s) dW(s) + dM^n(s)\right]\\
&&\quad - \int_0^t f_n(s,Y^n_s,Z^n_s) \phi(F(s,X_s,A_s)) ds + \int_0^t Y^n_{s} \phi^\prime(F_s(X_s,A_s)) h(s,X(s))A(s)  ds.
\end{eqnarray*}

\vspace{0.5cm}
Other examples are given by \cite[Examples 4 and 5]{cont:four:13}, namely
$$F(t,x_t,v_t) = x(t)^2 - \int_0^t v(u) du,\qquad F(t,x_t,v_t) = \exp \left( x(t) - \frac{1}{2} \int_0^t v(u) du\right).$$
Conditions on $b$ and $\sigma$ can be easily found such that (C3) and (C4) hold, especially if $X$ is given by \eqref{eq:Markov_SDE}.

\vspace{0.5cm}
Let us finish with the weak Euler-Maruyama scheme as in \cite{cont:lu:16}. We still consider the SDE \eqref{eq:SDE} with $b=0$ and the non-anticipative functional $X^n$ given by the recursion
$$X^n(t_{j+1}) = X^n(t_j) + \sigma(t_j,X^n_{t_j}) (W(t_{j+1}) - W(t_{j})).$$
For a Lipschitz functional $h : D([0,T],\bR^d) \to \bR$, consider the ``weak Euler approximation'' 
$$F_n(t) = \bE \left[ g(X^n_T)| \cF^W_t\right]$$
of the conditional expectation $\bE \left[ g(X_T)| \cF^W_t\right]$, where $\bF^W$ is the filtration generated by the Brownian motion $W$. This weak approximation is computed by initializing the scheme on $[0,t]$ with $\omega$ (a path of the Brownian motion) and then iterating the scheme with the increments of the Wiener process between $t$ and $T$. Then $F_n \in \bC^{1,\infty}_{{\rm loc}}$ (see \cite[Theorem 3.1]{cont:lu:16}). Moreover since we have a martingale, $\Theta_1(s)=0$. Under our setting and thanks to \cite[Theorem 4.1]{cont:lu:16}, (C3) holds. (C4) does not hold on the whole interval $[0,T]$. Nevertheless this functional is locally regular (\cite[Definition 7]{cont:lu:16}) and on our neighbourhood of $T$, one can easily get (C4) provided that $g$ is bounded for example.


\bibliographystyle{plain}
\bibliography{biblio_sing_BSDE_jumps}

\end{document}